\documentclass[12pt]{article}

\usepackage{amsthm,amsfonts,amsmath,amssymb,amscd}
\usepackage{bm}
\usepackage{mathabx}

\newcommand*{\QEDB}{\hfill\ensuremath{\blacksquare}}
\usepackage{authblk}

\begin{document}
	\title{Stability in the process with   lag of interaction between mining and processing industries}
	\author[1]{A.I.Kirjanen\thanks{kirjanen@mail.ru}}
	\author[1]{O.A.Malafeyev \thanks{malafeyevoa@mail.ru}}
	\author[1]{N.D.Redinskikh \thanks{redinskich@yandex.ru}}
	\affil[1]{Saint-Petersburg State University, Faculty of Applied Mathematics and Control Processes, Universitetskii prospekt 35, Petergof, Saint Petersburg, Russia, 198504}
	
\date{}
\maketitle

\begin{abstract}
A mathematical model of dynamic interaction between mining and processing industries is formalized and studied in the paper. The process of interaction is  described by a system of two delay differential equations. The criterion for asymptotic stability of nontrivial equilibrium point is obtained when both industries co-work steadily. The problem is reduced to finding stability criterion for  quasi-polynomial of second order. Time intervals between deliveries of raw materials which make it possible to preserve stable interaction between the two industries are found.
\end{abstract}

\section{Introduction}

A mathematical model of dynamic interaction between mining and processing industries is formalized and studied in the paper. It is supposed, that the resource is mined by the first  industry and then it is  transformed into some product by the second industry.  Let us denote the amount of the resources mined by $P$  and the number of companies, producing the final product from this resource. Similarly, it is possible to consider $Q$ as the total output in the processing industry.

Let us assume, that the volume of the mined resources  is growing with a coefficient $a>0$ due to ongoing mineral exploration and the amount of raw materials for processing industry is unlimited. We assume that the increase in the number of processing companies leads to the reduction of the volume of extracted raw materials with coefficient $b>0$ and, conversely, an increase in the  amount of raw materials implies an increase in the number of processing companies with coefficient $d>0$. Although the extraction of raw materials occurs continuously, it is shipped to producers in portions with some positive time lag h. So the amount of raw materials mined depends on its volume obtained earlier at time moment $(t-h)$ when the last supply was made. If this quantity of ``old" mined materials is large, the rate of its extraction must decrease with coefficient $e>0$.

In the absence of raw materials the number of processing companies is reduced with the coefficient $c>0$. In the model described by system (\ref{EQ1}) the number of processing companies at time moment $t$, also depends on the number of already operating companies at time moment $(t-h)$.

So we get that the process of mutual interaction between industries can be described by the differential equation system as follows
\begin{eqnarray}\label{EQ1}\frac{dP(t)}{dt}=(a-eP(t-h)-bQ(t))P(t) \nonumber\\
\frac{dQ(t)}{dt} =(-c+dP(t)-fQ(t-h))Q(t) . \end{eqnarray}
This system has a nontrivial equilibrium
\begin{equation}\label{EQ2}P^{*} =\frac{af+bc}{bd+ef} ;{\rm \; \; Q}^{*} =\frac{ad-ce}{bd+ef},\end{equation}
provided that \begin{equation}\label{EQ3}ad>ce.\end{equation}
Our aim is to find conditions on the coefficients of system (\ref{EQ1}) for stabilizing the equilibrium point (\ref{EQ2}). This would mean that the first industry produces such a quantity of raw materials that they will be processed by the second industry. In this case, it will not be overstocking in the warehouses and there will be sufficient volumes of raw materials.

Necessary and sufficient criterion for stable coexistence of two competitors was obtained in [1]. Graphically the stability areas were described  in the form of multidimensional cones in [2]. In this article, our aim is to express conditions for asymptotic stability in the form of inequalities. In this case these conditions help to solve problems of control and stabilization. A delay effect on the stability of the equilibrium point was studied in [3,4,5].

A number of mathematical models describing the interaction between agents based on the game theory was considered in [6--18]. The results on the business security, the impact of external factors on the growth of the business are given in [19,20]. Mathematical models  for delay-dependent linear systems with multiple time delays, for growing tumor, exponential stability with several delays, stability criteria for high even order delay differential equations were considered in [21--52].

\section{ Stability criterion for  quasi-polynomial}
By changing variables $x=P-P^{*} ,{\rm \; \; \; y}=Q-Q^{*} $ in the system (\ref{EQ1}) and writing down a linear approximation system, we get the characteristic quasi-polynomial and the characteristic equation as follows
\begin{equation}\label{EQ4}H(z)=z^{2} e^{2z} +c_{1} e^{2z} +c_{2} ze^{z} +c_{3} =0.\end{equation}
Here $c_{1} =bdh^{2} P^{*} Q^{*} ,{\rm \; \; }c_{2} =(eP^{*} +fQ^{*} )h,{\rm \; \; }c_{3} =efh^{2} P^{*} Q^{*}$. To get the conditions under which the roots of the quasi-multinomial (\ref{EQ4}) lie in the left half-plane, we use Pontryagin and Hermite - Biehler criteria [1,3,4,27,28].

\textbf{Theorem 1.} The roots of the quasi-polynomial (\ref{EQ4}) with positive coefficients lie in the left half-plane if one of two following assertions $A$ or $B$ is fulfilled:
	
	Assertion A:
	
	I.1. $0<c_{1} <\pi ^{2} $;
	
	I.2. $0<c_{2} <\frac{z(\dot{y})}{\dot{y}} =\frac{2(\dot{y}^{2} -c_{1} )\sin \dot{y}}{\dot{y}} $, here $\dot{y}\in (\sqrt{c_{1} } ,\pi )$ is a unique root of the equation $\tan y=\frac{(c_{1} -y^{2} )y}{c_{1} +y^{2} } $; notice, that   $\dot{y}\in (\frac{\pi }{2} ;\pi );$
	
	I.3.1. If $y_{1} \in (0;\frac{\pi }{2} )$ is a unique root of the equation $2(y^{2} -c_{1} )\sin y=c_{2} y$ then the following conditions are fulfilled:
	
	I.3.1.1. $c_{1} +c_{3} <y_{1} ^{2} $ and   I.3.1.2. $c_{1} +c_{2} \frac{\pi }{2} <c_{3} +\frac{\pi ^{2} }{4} $.
	
	I.3.2. If${\rm \; \; \; \; y}_{{\rm 1}} ,y_{2} \in (\frac{\pi }{2} ;\pi )$, $y_{1} <y_{2} $ are the roots of the equation $2(y^{2} -c_{1} )\sin y=c_{2} y$ then the following  conditions are fulfilled:
	
	I.3.2.1.  $y_{1} ^{2} <c_{1} +c_{3} <y_{2} ^{2} $ and   I.3.2.2.  $c_{1} +c_{2} \frac{\pi }{2} >c_{3} +\frac{\pi ^{2} }{4} $
	
	Assertion B:

	II.1. $\pi ^{2} <c_{1} <4\pi ^{2} $;
	
	II.2. $0<c_{2} <\frac{z(\dot{y})}{\dot{y}} =\frac{2(\dot{y}^{2} -c_{1} )\sin \dot{y}}{\dot{y}} $, here $\dot{y}\in (\pi ;\sqrt{c_{1} } )$ is the unique root of the equation $\tan y=\frac{(c_{1} -y^{2} )y}{c_{1} +y^{2} } $; notice, that   $\dot{y}\in (\pi ;\frac{3\pi }{2} );$
	
	II.3.1. If $y_{1} <y_{2} $ are the roots of the equation $2(y^{2} -c_{1} )\sin y=c_{2} y$ and$y_{1} ,C_{2} \in (\pi ;\sqrt{c_{1} } )$, then the following  conditions are fulfilled:
	
	II.3.1.1.$y_{1} ^{2} <c_{1} +c_{3} <y_{2} ^{2} $ and   II.3.1.2.$c_{1} +c_{2} \frac{\pi }{2} >c_{3} +\frac{\pi ^{2} }{4} $.
	
	II.3.2. If $y_{2} \in (\frac{3\pi }{2} ;\sqrt{c_{1} } )$ is the unique root of the equation $2(y^{2} -c_{1} )\sin y=c_{2} y$ then the following conditions are fulfilled:
	
	II.3.2.1. $c_{1} +c_{3} >y_{2} ^{2} $ and   II.3.2.2. $\frac{9\pi ^{2} }{4} +c_{2} \frac{3\pi }{2} +c_{3} <c_{1} $.

\begin{proof} Let us substitute $z=iy$ into the quasi-polynomial (\ref{EQ4}) and write down its real $F(y)$ and imaginary $G(y)$ parts:
	\begin{equation}\label{EQ5}F(y)=(-y^{2} +c_{1} )\cos 2y-c_{2} y\sin y+c_{3} , \end{equation}
	\begin{equation}\label{EQ6}G(y)=(-2y^{2} \sin y+2c_{1} \sin y+c_{2} y)\cos y=(c_{2} y-z(y))\cos y=0.\end{equation}
	Then derivative $G'(y)$ may be written as follows
	\begin{equation}\label{EQ7}G'(y)=(c_{2} -z'(y))\cos y-(c_{2} y-z(y))\sin y.\end{equation}
	Here $z(y)=2(y^{2} -c_{1} )\sin y$.
	
	From Pontryagin and Hermite - Biehler criteria it is known that the roots of quasi-polinomial (\ref{EQ4}) have negative real parts if the vector of gain-phase cha-racteristic (amplitude-phase characteristic) $w=H(iy)$ monotonically rotates counterclockwise round the origin with positive rate [1,3,4,27,28]. It means that the gain-phase characteristic turning around the origin crosses every line passing through the point (0,0) at the positive angle without touching it [1,3]. In this case all roots of the functions $F(y)$ and $G(y)$ are real, simple, alternate and the inequality
	\begin{equation}\label{EQ8} F(y)G'(y)-F'(y)G(y)>0 \end{equation}
	holds for all $y$. For the stability of quasi-polinomial (\ref{EQ4}) it is sufficient for the inequality (\ref{EQ8}) to be valid only at the roots of function $G(y)$. In the future, we will consider the inequality
	\begin{equation} \label{EQ9} F(y)G'(y)>0 \end{equation}
	at the roots of the function $G(y)$.
	
	Due to the Pontryagin criterion [28] inequalities (\ref{EQ8}) and (\ref{EQ9}) are valid iff the rotation angle of gain-phase characteristic $(F(y);G(y))$ around the origin asymptotically tends to $\varphi (-2\pi k+\varepsilon \le y\le 2\pi k+\varepsilon )\approx (4ks+r)\pi $ as $k$ tends to infinity. Here $s$ is the degree of quasi-polynomial (\ref{EQ4}) with respect to $e^{z} $and $r$ is the degree of quasi-polynomial (\ref{EQ4}) with respect to $z$. In our case $s=r=2$ so, polynomials $F(y)$ and $G(y)$ have each 10 roots on the segment $[-2\pi+\varepsilon ;2\pi +\varepsilon ]$. As the function $G(y)$ is odd it can  not have more than 9 roots over the segment $[-2\pi ;2\pi ]$. One root is $y_{0} =0$, the other ones should lie symmetrically with 4 on each side. The tenth root will be discussed later.
	
	\textbf{Proof of Assertion A.} Note that $\frac{\pi }{2}$ and $\frac{3\pi }{2}$ are the roots of the function $G(y)$. Other two roots $y_{1} <y_{2} \in (0;2\pi )$ are the solutions of the equation
	\begin{equation}\label{EQ10}z(y)=2(y^{2} -c_{1} )\sin y=c_{2} y.\end{equation}
	Since $z'(0)=-2c_{1} <0$, the equation (\ref{EQ10}) has no roots for small positive y. On the other hand, the function $z(y)$ crosses the x-axis at  the points $y=\sqrt{c_{1} }$ and $y=\pi$, so two cases can occur.
	
	Let $0<c_{1} <\pi ^{2} $. Then the function $z(y)=2(y^{2} -c_{1} )\sin y$ is positive over $(\sqrt{c_{1} } ,\pi )$ and equation (\ref{EQ10}) has 2 solutions on this interval if the coefficient $A_{2}$ is less than the slope of tangent line to the graph of the function $z=z(y)$, drawn from the origin. Then the touch point $\dot{y}$ is the solution of the equation $z(y)=z'(y)y$ or $\tan y=\frac{(c_{1} -y^{2} )y}{c_{1} +y^{2} } $. As we consider the solution of this equation on $(\sqrt{c_{1} } ,\pi )$, so $\tan y<0$ and the touch point $\dot{y}\in (\frac{\pi }{2} ;\pi )$. In this case equation (\ref{EQ10})  has two roots $y_{1} <y_{2} \in $$(0;\pi )$ if the condition I.2 of the Theorem 1 is satisfied.
	
	To satisfy the Hermite - Biehler conditions  it is necessary to have alternation of signs both of derivative $G'(y)$ and function  $F(y)$ in the roots of $G(y)$ and their multiplication should satisfy condition (\ref{EQ9}).
	
	For the  root $y_0=0$  we get $G'(0)=2c_{1} +c_{2} >0,F(0)=c_{1} +c_{3} >0$. Let's assume that other roots of the function $G(y)$ are ordered as follows: $y_{1} <\frac{\pi }{2} <y_{2} <\frac{3\pi }{2} $. Then $c_{2} <z'(y_{1} )$ and the inequality (\ref{EQ7}) is equivalent to the condition  $G'(y_{1} )<0$. For the root $y_{2}$ the inverse inequality         is  valid. But $\cos y_{2} <0$ so  $G'(y_{2} )<0$. Between these two roots a straight line $z=c_{2} y$ lies under the graph of the function $z=z(y)$ so the inequality $c_{2} y<z(y)$  is valid at the point $y=\frac{\pi }{2}$. From (\ref{EQ7}) we get   $G'(\frac{\pi }{2} )>0$. Taking into account the condition I.1, we obtain   $G'(\frac{3\pi }{2} )>0$. So the signs of the derivative $G'(y)$ in roots of function $G(y)$ alternate.    The function $F(y)$ has the form $F(y)=c_{1} +c_{3} -y_{j} ^{2} ,{\rm \; j}=1,2$  at the roots $y_{1} $ and $y_{2} $ of   the function $G(y)$. If $c_{1} +c_{3} <y_{1} ^{2} $ then the both inequalities $F(y_{j} )<0,{\rm \; \; }j=1,2$ are valid. Inequality $F(\frac{\pi }{2} )>0$ corresponds to the condition I.3.1.2 of the theorem 1 and then inequality $F(\frac{3\pi }{2} )>0$ is valid as well. So the signs of the function $F(y)$ alternate in roots of $G(y)$.
	
	We can get the tenth root of $G(y)$ by shifting the segment $[-2\pi ;2\pi ]$ to the right so that the root of the function $G(y)$ $y_{3} \in (2\pi ;2\pi +\frac{\pi }{2} )$ will be in this segment. As there are no roots of the function $G(y)$ on interval $(-2\pi ;-\frac{3\pi }{2} )$ we don't lose any of the root of the function $G(y)$ on the left side. So we have constructed a segment of $4\pi$-length, and there are 10 roots of the function $G(y)$ on this segment.
	
	Let the roots of the function $G(y)$ be ordered as follows:  $\frac{\pi }{2} <y_{1} <y_{2} <\frac{3\pi }{2} $. In this case a straight line$z=c_{2} y$ is located over graph of the function $z=z(y)$ on the interval $(0;y_{1} )$, so at point $y=\frac{\pi }{2} $   inequality $c_{2} \frac{\pi }{2} >z(\frac{\pi }{2} )$ is valid. From (\ref{EQ7}) we get $G'(\frac{\pi }{2} )<0$. Similarly to the previous case, $c_{2} <z'(y_{1} )$ and taking into account the inequality $\cos y_{1} <0$ from (\ref{EQ7}) we obtain $G'(y_{1} )>0$. The inverse inequality $c_{2} >z'(y_{2} )$ with $\cos y_{2} <0$ yields $G'(y_{2} )<0$. Inequality $G'(\frac{3\pi }{2} )>0$ is also valid, and we get the alternation of signs of the derivative $G'(y)$ in the roots of the function $G(y)$. From the condition I.3.2.2 of the Theorem 1 we obtain $F(\frac{\pi }{2} )<0$   and from inequalities $y_{1} ^{2} <c_{1} +c_{3} <y_{2} ^{2} $ we obtain $F(y_{1})>0$ and $F(y_{2} )<0$. Inequality $F(\frac{3\pi }{2} )>0$ is valid as well. So the signs of the function $F(y)$ in roots of $G(y)$ alternate. Further, all $2\pi$-long segments will include four roots of the function $G(y)$ and there will be a similar alternation of the signs of the derivative  $G'(y)$ and the signs of the function  $F(y)$ in these roots  and  inequality (\ref{EQ9})  is valid as well. Assertion A is proved.

	\textbf{Proof of Assertion B.} Let $\pi ^{2} <c_{1} <4\pi ^{2} $. Then function $z(y)=2(y^{2} -c_{1} )\sin y$ is positive on $(\pi ;\sqrt{c_{1} } )$ and equation (\ref{EQ10}) has 2 solutions on this interval if coefficient $c_{2} $ is less than the slope of tangent to function graph $z=z(y)$, drawn from the origin. It was noted at the proof of the assertion À that the touch point $\dot{y}$ is the solution of the equation $\tan y=\frac{(c_{1} -y^{2} )y}{c_{1} +y^{2} } $. We consider the solution of this equation on $(\pi ;\sqrt{c_{1} } )$, so the touch point $\dot{y}\in (\pi ;\frac{3\pi }{2} )$.
	
	Let the roots of $G(y)$ are ordered as follows   $\frac{\pi }{2} <y_{1} <y_{2} <\frac{3\pi }{2} $. We have considered such sequence of roots in the previous part of the proof but now $\pi <y_{1} $. As before we have $c_{2} \frac{\pi }{2} >z(\frac{\pi }{2} )$ and from (\ref{EQ7})  we get $G'(\frac{\pi }{2} )<0$. At the point $y_{1} $  we have $c_{2} <z'(y_{1} )$ and $\cos y_{1} <0$ so the inequality (\ref{EQ7}) is equivalent to the inequality $G'(y_{1} )>0$. At the point $y_{2} $  we have $c_{2} >z'(y_{2} )$ and $\cos y_{2} <0$, so from (\ref{EQ7}) we obtain $G'(y_{2} )<0$. Inequality $G'(\frac{3\pi }{2} )>0$ is fulfilled as well. Inequalities  $F(\frac{\pi }{2} )<0$, $F(y_{1} )>0$, $F(y_{2} )<0$, $F(\frac{3\pi }{2} )>0$  correspond to the inequalities  $c_{1} +c_{2} \frac{\pi }{2} <c_{3} +\frac{\pi ^{2} }{4} $, $y_{1} ^{2} <c_{1} +c_{3} <y_{2} ^{2} $, so the signs of the function $F(y)$ in roots of $G(y)$ alternate. If $\frac{\pi }{2} <y_{1} <\frac{3\pi }{2} <y_{2} $, and $\pi <y_{1} $, then $G'(\frac{\pi }{2} )<0$. The root $y=\frac{3\pi }{2} $ lies between two roots $y_{1} <y_{2} $ so the straight line $z=c_{2} y$ is located under the graph of the function $z=z(y)$ or $c_{2} y<z(y)$. From this inequality we obtain $G'(y)=-\sin y{\rm \; }(c_{2} y-z(y)<0$ for $y=\frac{3\pi }{2} $. So the signs of the derivative $G'(y)$ in roots of function $G(y)$ alternate.  The corresponding alternation of signs of function $F(y)$ takes place iff  the inequalities $\frac{9\pi ^{2} }{4} +c_{2} \frac{3\pi }{2} +c_{3} <c_{1} $  and $c_{1} +c_{3} >y_{2} $ are satisfied. If $\frac{\pi }{2} <\frac{3\pi }{2} <y_{1} <y_{2} $, then the both inequalities $c_{1} +c_{3} <y_{1} $ and $c_{1} +c_{3} >y_{2} $ are fulfilled together. It contradicts  to the inequality $y_{1} <y_{2} $. So  the assertion B of the  Theorem 1 is proved.
	
	If $c_{1} >4\pi ^{2}$, the segment $[0;2\pi]$ contains only two roots of the function  $G(y)$, whereas Pontryagin criterion  requires four roots.\QEDB
\end{proof}
\section{Criterion for stable co-functioning of the two industries}
In this paragraph we obtain conditions under which the goods quantity meets the demand. Under these conditions there is as overproducing of raw materials and its shortage as well. Proving  the theorem 1 we have considered only linear system associated with the equilibrium point  $(P^{*} ;Q^{*})$ without square members. From the theory of the differential equations it is known that in this case the equilibrium point of system (\ref{EQ1}) is asymptotically stable if the approximating linear system is asymptotically  stable. The last statement is valid if the roots of the quasi-polynomial (\ref{EQ4}) lie in the left half-plane.    Thus, we obtain the following theorem.

\textbf{Theorem 2.} A non-trivial equilibrium $(P^{*} ;Q^{*} )$ of the system (\ref{EQ1}) is asymptotically stable iff the following conditions are fulfilled:
	
	Assertion A:
	
	I.1. $0<bdh^{2} P^{*} Q^{*} <\pi ^{2} $
	
	I.2. $0<(eP^{*} +fQ^{*} )h<\frac{z(\dot{y})}{\dot{y}} =\frac{2(\dot{y}^{2} -bdh^{2} P^{*} Q^{*} )\sin \dot{y}}{\dot{y}} $, here $\dot{y}\in (\sqrt{bdh^{2} P^{*} Q^{*} } ,\pi )$ is the unique root of the equation $\tan y=\frac{(bdh^{2} P^{*} Q^{*} -y^{2} )y}{bdh^{2} P^{*} Q^{*} +y^{2} } $;
	
	I.3.1. If $y_{1} \in (0;\frac{\pi }{2} )$ is a root of equation $2(y^{2} -bdh^{2} P^{*} Q^{*} )\sin y=(eP^{*} +fQ^{*} )hy$ then the following conditions are fulfilled:
	
	I.3.1.1. $bdh^{2} P^{*} Q^{*} +efh^{2} P^{*} Q^{*} <y_{1} ^{2} $ and
	
	I.3.1.2. $bdh^{2} P^{*} Q^{*} +(eP^{*} +fQ^{*} )h\frac{\pi }{2} <efh^{2} P^{*} Q^{*} +\frac{\pi ^{2} }{4} $.
	
	I.3.2.  If $y_{1} ,y_{2} \in (\frac{\pi }{2} ;\pi )$, $y_{1} <y_{2} $ are the roots of the equation $2(y^{2} -bdh^{2} P^{*} Q^{*} )\sin y=(eP^{*} +fQ^{*} )hy$ then the following  conditions are fulfilled:
	
	I.3.2.1. $y_{1} ^{2} <bdh^{2} P^{*} Q^{*} +efh^{2} P^{*} Q^{*} <y_{2} ^{2} $ and
	
	I.3.2.2. $bdh^{2} P^{*} Q^{*} +(eP^{*} +fQ^{*} )h\frac{\pi }{2} >efh^{2} P^{*} Q^{*} +\frac{\pi ^{2} }{4} $
	
	Assertion B:
	
	A non-trivial equilibrium $(P^{*} ;Q^{*})$ of the system (\ref{EQ1}) is asymptotically stable if the following conditions are fulfilled:
	
	II.1. $\pi ^{2} <bdh^{2} P^{*} Q^{*} <4\pi ^{2} $;
	
	II.2. $0<(eP^{*} +fQ^{*} )h<\frac{z(\dot{y})}{\dot{y}} =\frac{2(\dot{y}^{2} -bdh^{2} P^{*} Q^{*} )\sin \dot{y}}{\dot{y}} $,here $\dot{y}\in (\pi ;\sqrt{bdh^{2} P^{*} Q^{*} } )$ is the unique root of the equation $\tan y=\frac{(bdh^{2} P^{*} Q^{*} -y^{2} )y}{bdh^{2} P^{*} Q^{*} +y^{2} } $;
	
	II.3.1. If $y_{1} ,y_{2} \in (\pi ;\frac{3\pi }{2} )$ $y_{1} <y_{2} $ are the roots of the equation $2(y^{2} -bdh^{2} P^{*} Q^{*} )\sin y=(eP^{*} +fQ^{*} )hy$  then  the following  conditions are valid:
	
	II.3.1.1. $y_{1} ^{2} <bdh^{2} P^{*} Q^{*} +efh^{2} P^{*} Q^{*} <y_{2} ^{2} $ and
	
	II.3.1.2. $bdh^{2} P^{*} Q^{*} +(eP^{*} +fQ^{*} )h\frac{\pi }{2} >efh^{2} P^{*} Q^{*} +\frac{\pi ^{2} }{4} $.
	
	II.3.2.  If $y_{2} \in (\frac{3\pi }{2} ;\sqrt{bdh^{2} P^{*} Q^{*} } )$ is a root of the equation$2(y^{2} -bdh^{2} P^{*} Q^{*} )\sin y=(eP^{*} +fQ^{*} )hy$ then the following conditions are valid:
	
	II.3.2.1. $bdh^{2} P^{*} Q^{*} +efh^{2} P^{*} Q^{*} >y_{2} ^{2} $ ,
	
	II.3.2.2. $\frac{9\pi ^{2} }{4} +(eP^{*} +fQ^{*} )h\frac{3\pi }{2} +efh^{2} P^{*} Q^{*} <bdh^{2} P^{*} Q^{*} $.

\section{Conclusions}

1. A mathematical model of dynamic interaction between mining and processing industries is described by the system of two nonlinear delay differential equations. At the proposed model we take into account  the volume of raw materials mined and shipped in the preceding time $(t-h)$. In the Theorem 2 we find conditions on the coefficients of system (\ref{EQ1}) for stabilizing the equilibrium point (\ref{EQ2}). It means that the first industry produces such a quantity of raw materials that they will immediately be processed by the second industry. It means that there is a balance  between the amount of extracted raw materials and the number of processing enterprises.

2. The conditions under which the solutions of system (\ref{EQ1}) tend to equilibrium point $(Q^{\bullet };P^{\bullet })$  are given in Theorem 2. The conditions I.1, I.2, I.3.1.1 and I.3.1.2 taken together, give a size of admissible interval between deliveries of raw material at which the balance between  industries is kept. The size of this $h^{2} <\frac{y_{1} ^{2} }{(bd+ef)P^{\bullet } Q^{\bullet } } $, where $y_{1} \in (0;\frac{\pi }{2} )$ is a root of equation $2(y^{2} -bdh^{2} P^{*} Q^{*} )\sin y=(eP^{*} +fQ^{*} )hy$.

3. If the interval between deliveries of raw material is not small but conditions I.1, I.2, I.3.2.1 and I.3.2.2 of Theorem 2 taken together are fulfilled then the balance between  industries is kept also. In this case size of this admissible interval is $\frac{y_{1} ^{2} }{(bd+ef)P^{\bullet } Q^{\bullet } } <h^{2} <\frac{y_{2} ^{2} }{(bd+ef)P^{\bullet } Q^{\bullet } } $, where $y_{1} ,y_{2} \in (\frac{\pi }{2} ;\pi )$, $y_{1} <y_{2} $    are the roots of the equation $2(y^{2} -bdh^{2} P^{*} Q^{*} )\sin y=(eP^{*} +fQ^{*} )hy$.

4. Moreover if delay $h$ is rather large $\pi ^{2} <bdh^{2} P^{*} Q^{*} <4\pi ^{2} $ (condition II.I of Theorem 2) the stability of the equilibrium point $(Q^{\bullet };P^{\bullet } )$ can be restored if conditions of Assertion B of Theorem 2 are fulfilled.

5. If this time interval between  supplies of raw materials increases, the asymptotic stability of equilibrium point stabilization. Restrictions obtained in the Theorem 2 allow to solve the problem of stabilization of equilibrium point.

6. Assume that a priori it is known the quantity of final products required to meet demand in the region, and to manifacture these final products we need a certain amount of raw materials. It means that initially the equilibrium point is known, but the point $(Q^{\bullet} ;P^{\bullet})$ determined by the formula (\ref{EQ2}), does not meet the demand for these products. In this case, it is necessary to change the technological process of extraction and processing (change coefficients $a,b,c,d,e,f$) so that point $(Q^{\bullet }; P^{\bullet })$ is consistent with the economically justified demand. Then it is possible to solve the above mentioned problem of finding time intervals between deliveries of raw materials which make it possible to preserve stable interaction between the two industries.

\end{document}